# AN EXTENDED NEWTON-TYPE METHOD IN ℝ AND POLYNOMIOGRAPHY VIA DIFFERENT ITERATIVE METHODS

Nazli Karaca and Isa Yildirim

**Abstract:** The aim of this paper is to introduce a new Newton-type iterative method and then to show that this process converges to the unique solution of the scalar nonlinear equation $f(x) = 0$ under weaker conditions involving only $f$ and $f'$ by fixed point techniques. Also, by using this iteration process quite new nicely looking graphics are obtained.

## 1. Introduction and Preliminaries

Newton's method or Newton-Raphson method, as it is generally called in the case of scalar equations $f(x) = 0$, is one of the most used iterative procedures for solving such nonlinear equations. Newton's method is defined by an iterative sequence

$$x_{n+1} = x_n - \frac{f(x_n)}{f'(x_n)}, n \geq 0, \tag{1}$$

under suitable assumptions on $f$ and $f'$. Note that (1) can be viewed as the sequence of successive approximations (Picard iteration) of the Newton iteration function given by

$$T(x) = x - \frac{f(x)}{f'(x)}$$

where $f' \neq 0$. Moreover, under appropriate conditions, $x^*$ is a solution of $f(x) = 0$ if and only if $x^*$ is a fixed point of the iteration function $T$.

There exist several convergence theorems under weak conditions which involve $f$, $f'$ and $f''$ in literature for the Newton's method, see for example [7], [11] and [18].

**Theorem 1.** *([7]) Let $f: [a, b] \to \mathbb{R}$. Suppose that the following conditions hold:*

*(i) $f(a)f(b) < 0$;*

*(ii) $f \in C^2[a, b]$ and $f'(x)f''(x) \neq 0$, $x \in [a, b]$;*

*Then the sequence $\{x_n\}$ defined by (1) starting with an initial guess $x_0 \in [a, b]$ converges to $x^*$; the unique solution of $f(x) = 0$ in $[a, b]$. Moreover, we have the following estimation*

$$|x_n - x^*| \leq \frac{M_2}{2m_1} |x_n - x_{n-1}|, n \geq 1, \tag{2}$$

*holds, where*

$$m_1 = \min_{x \in a,b]} |f'(x)| \text{ and } M_2 = \max_{x \in a,b]} |f''(x)|.$$

For numerical point of view, Theorem 1 is widely applicable but there exist more general results based on weaker smoothness conditions. In a series of papers [2]-[4], Berinde obtained more general convergence results which extend Newton's method both scalar ([2], [3]) and $n$-dimensional equations [4]. These results can be applied to weakly smooth functions. The term extended Newton method was adopted in view of the fact that the iterative process (1) has been extended from $[a, b]$ to the whole real axis $\mathbb{R}$.

One of the scalar variant of these results is stated below.

**Theorem 2.** *([2]) Let $f: [a, b] \to \mathbb{R}$, where $a < b$. If the following conditions hold:*

*($f_1$) $f(a)f(b) < 0$;*

*($f_2$) $f \in C^1[a, b]$ and $f'(x) \neq 0$, $x \in [a, b]$;*

*($f_3$) $2m > M$, where*

$$m = \min_{x \in [a,b]} |f'(x)| \text{ and } M = \max_{x \in [a,b]} |f'(x)|. \tag{3}$$

*Then the Newton iteration $\{x_n\}$, defined by (1) starting with $x_0 \in [a,b]$ converges to $x^*$; the unique solution of $f(x) = 0$ in $[a,b]$. Moreover, the following estimation*

$$|x_n - x^*| \leq \frac{M}{m} |x_n - x_{n+1}|, n \geq 0, \tag{4}$$

*holds.*

All the proofs in [2], [3] are based on a classical technique which focuses on the behavior of the sequence $\{x_n\}$ defined in (1).

Recently, Sen et al. [20] extended Theorem 2 to the case of a Newton-like iteration of the form given as:

$$x_{n+1} = x_n - \frac{2f(x_n)}{f'(x_n) + M_1 f(x_n)}, n \geq 0, \tag{5}$$

with $M_1 f(x) = sgn f'(x) \cdot M$, where $M$ is defined by (3).

Later, this result was extended to the $n$-dimensional case [21]. However, in both cases an extended Newton-like algorithm was used.

There exists a strong link of Newton's methods with iteration processes in fixed point theory. In 2007, Agarwal, O'Regan and Sahu [1] have introduced the S-iteration process as follows: Let $X$ be a normed space, $C$ a nonempty convex subset of $X$ and $T: D \to C$ an operator. Then, for arbitrary $x_0 \in C$, the S-iteration process is defined by

$$\begin{cases} x_{n+1} = (1 - \alpha_n) T x_n + \alpha_n T y_n, \\ y_n = (1 - \beta_n) x_n + \beta_n T x_n, n \geq 0, \end{cases} \tag{6}$$

where $\{\alpha_n\}$ and $\{\beta_n\}$ are sequences in $(0,1)$.

In 2009, Yildirim and Ozdemir [23] proved some convergence result by using the following iteration process: For an arbitrary fixed order $r \geq 2$,

$$\begin{cases} x_{n+1} = (1 - \alpha_{1n})y_{n+r-2} + \alpha_{1n}Ty_{n+r-2}, \\ y_{n+r-2} = (1 - \alpha_{2n})y_{n+r-3} + \alpha_{2n}Ty_{n+r-3}, \\ \vdots \\ y_{n+1} = (1 - \alpha_{(r-1)n})y_n + \alpha_{(r-1)n}Ty_n, \\ y_n = (1 - \alpha_{rn})x_n + \alpha_{rn}Tx_n, \quad n \geq 0, \end{cases} \quad (7)$$

or, in short,

$$\begin{cases} x_{n+1} = (1 - \alpha_{1n})y_{n+r-2} + \alpha_{1n}Ty_{n+r-2}, \\ y_{n+r-i} = (1 - \alpha_{in})y_{n+r-(i+1)} + \alpha_{in}Ty_{n+r-(i+1)}, \\ y_n = (1 - \alpha_{rn})x_n + \alpha_{rn}Tx_n, \quad n \geq 0, \end{cases} \quad (8)$$

where $\{\alpha_{1n}\}$ and $\{\alpha_{in}\}$, $i = 2, \ldots r$, are real sequence in $[0,1)$.

**Remark 1.** *i) If we take $r = 2$ in (8), we obtain the iteration process in [22].*

*ii) If we take $r = 3$ in (8), we obtain the following iteration process:*

$$\begin{cases} x_{n+1} = (1 - \alpha_n)z_n + \alpha_n T(z_n), \\ z_n = (1 - \beta_n)y_n + \beta_n T(y_n), \\ y_n = (1 - \gamma_n)x_n + \gamma_n T(x_n), n \geq 0. \end{cases} \quad (9)$$

After, Khan [14] introduced a new iteration process for nonexpansive mappings, which he called 'Picard-Mann hybrid iteration process' and the convergence process is faster than Picard and Mann iteration process. Let $X$ be a normed space, $C$ a nonempty convex subset of $X$ and $T: C \to C$ an operator. Then, for arbitrary $x_0 \in C$, the Picard-Mann hybrid iteration process is defined by

$$\begin{cases} x_{n+1} = Ty_n, \\ y_n = (1 - \alpha_n)x_n + \alpha_n Tx_n, n \geq 0, \end{cases} \quad (10)$$

where $\{\alpha_n\} \subset (0,1)$.

Karaca *et al* [12] obtained a convergence result for the iteration process (10) of Newton-like and they showed this iteration process is better than the Newton method (1) and the extended Newton-like method (5).

Recently, Kadioglu and Yildirim introduced an iteration process in [13]:

$$\begin{cases} x_{n+1} = Ty_n, \\ y_n = (1 - \alpha_n)z_n + \alpha_n Tz_n, \\ z_n = (1 - \beta_n)x_n + \beta_n Tx_n, n \geq 0, \end{cases} \quad (11)$$

where $\{\alpha_n\}$ and $\{\beta_n\}$ are sequences in $(0,1)$. And, they showed that the iteration process (11), for contractions, is faster than both the S-iteration process and the Picard-Mann hybrid iteration process.

Motivated by Newton's method and the other iteration process, we will introduce the iteration process (11) of Newton-like for a real-valued function $f$ defined on an open interval $I$ as follows: For arbitrary $x_0 \in I$, the iteration process (11) of Newton-like is defined by

$$\begin{cases} x_{n+1} = N(y_n), \\ y_n = (1 - \alpha_n)z_n + \alpha_n N(z_n), \\ z_n = (1 - \beta_n)x_n + \beta_n N(x_n), n \geq 0, \end{cases} \quad (12)$$

where $\alpha, \beta \in (0,1)$ and $N(x) = x - \frac{f(x)}{f'(x)}$.

The purpose of this paper is to prove that the iteration process (12) converges to the unique solution of the scalar nonlinear equation $f(x) = 0$ under weaker conditions involving only $f$ and $f'$. Also, by using this algorithm quite new nicely looking polynomiographs are obtained.

The following definitions and lemma will be needed in the sequel.

**Definition 1.** *Let $(X, d)$ be a metric space. A mapping $T: X \to X$ is said to be*

*(i) contraction if there exists a constant $\lambda \in 0,1)$ such that for any $x$, $y \in X$, the following condition hold:*

$$d(T(x), T(y)) \leq \lambda d(x, y).$$

*(ii) quasi-contraction [19] if there exist a constant $\lambda \in 0,1)$ such that for any $x \in X$ and $x^* \in F(T)$, we have*

$$d(T(x), x^*) \leq \lambda d(x, x^*), \quad (13)$$

*where, $F(T) = \{x \in X : Tx = x\} \neq \emptyset$.*

The following lemma will be used in the proof of the main result of this paper.

**Lemma 1.** *[5] Let $(X, d)$ be a complete metric space and $T: X \to X$ a quasi-contractive operator with $x^* \in F(T)$. Then $x^*$ is the unique fixed point of $T$ and the Picard iteration $\{T^n(x_0)\}$ converges to $x^*$ for each $x_0 \in X$.*

## 2. Main Result

We start with the our main result.

**Theorem 3.** *Let $f: [a, b] \to \mathbb{R}$ be a function such that the following conditions are satisfied*

$(f_1)$ $f(a)f(b) < 0$;

$(f_2)$ $f \in C^1[a, b]$ and $f'(x) \neq 0$, $x \in [a, b]$;

$(f_3)$ $2m > M$, where

$$m = \min_{x \in [a,b]} |f'(x)| \text{ and } M = \max_{x \in [a,b]} |f'(x)|.$$

*Then*

*(i) The iteration process (12) starting with an arbitrary point $x_0$ in $[a, b]$ converges the unique solution $x^*$ of $f(x) = 0$ in $[a, b]$.*

*(ii) We have the following error estimate*

$$|x_n - x^*| \leq \frac{mM}{(1+\alpha)M^2 - \beta m^2} |x_n - x_{n+1}|, \qquad (14)$$

*for $n \geq 0$.*

**Proof.** (i) By conditions $(f_1)$ and $(f_2)$ it follows that the equation $f(x) = 0$ has a unique solution $x^*$ in $(a, b)$.

Suppose that $T: [a, b] \to \mathbb{R}$ is the Newton-like iteration function associated with $f$, that is

$T(x) = V U_\alpha K_\beta(x)$, where $V, U_\alpha, K_\beta: [a,b] \to \mathbb{R}$ are defined as:

$$V(x) = x - \frac{f(x)}{f'(x)}, \tag{15}$$

$$U_\alpha(x) = x - \alpha \frac{f(x)}{f'(x)},$$

$$K_\beta(x) = x - \beta \frac{f(x)}{f'(x)}.$$

Note that, $x^*$ is a solution of $f(x) = 0$ if and only if $x^*$ is a fixed point of $V$, $U_\alpha$ and $K_\beta$, that is

$$V(x^*) = U_\alpha(x^*) = K_\beta(x^*) = x^*.$$

From (15), we get

$$K_\beta(x) - x^* = x - \beta \frac{f(x)}{f'(x)} - x^* = x - x^* - \beta \frac{f(x)}{f'(x)}. \tag{16}$$

Since $x^*$ is a solution of $f(x)$,

$$f(x) = f(x) - 0 = f(x) - f(x^*).$$

Using condition ($f_2$) and the mean value theorem, we have

$$f(x) = f'(\bar{y})(x - x^*), \tag{17}$$

where $\bar{y} = x^* + \lambda(x - x^*)$, $0 < \lambda < 1$. From (16) and (17), we obtain

$$K_\beta(x) - x^* = (x - x^*)\left(1 - \beta \frac{f'(\bar{y})}{f'(x)}\right) \tag{18}$$

for all $x \in [a, b]$.

By condition ($f_2$), $f'$ preserves sign on $[a,b]$. That is, $f'(\bar{y})/f'(x) > 0$. Thus for any $x \in [a, b]$ and $\bar{y}$ between $x^*$ and $x$,

$$1 - K_\beta \frac{f'(\bar{y})}{f'(x)} < 1. \tag{19}$$

Also, using condition (f$_2$),

$$\beta \frac{f'(\bar{y})}{f'(x)} < \frac{f'(\bar{y})}{f'(x)} = \left|\frac{f'(\bar{y})}{f'(x)}\right| = \frac{|f'(\bar{y})|}{|f'(x)|} \leq \frac{M}{m} < 2,$$

which implies that

$$1 - \beta \frac{f'(\bar{y})}{f'(x)} > -1 \tag{20}$$

where, $x^* < \bar{y} < x$ and for all $x \in [a, b]$. From (19), (20) and the continuity of $f'$, we have

$$\lambda = \max_{x, \bar{y} \in [a,b]} \left|1 - \beta \frac{f'(\bar{y})}{f'(x)}\right| < 1 \text{ and } 0 < \lambda < 1,$$

which together with (18) implies that

$$|K_\beta(x) - x^*| \leq \lambda |x - x^*|, \ \forall x \in [a, b].$$

In a similar way we obtain

$$|U_\alpha(x) - x^*| \leq \delta |x - x^*|$$

and

$$|V(x) - x^*| \leq \gamma |x - x^*|$$

where $0 < \delta, \gamma < 1$. If we can use same arguments as given in the proof of Theorem 6 in [5] and obtain the following

$$x_{n+1} = K_\beta(x_n) \in [a, b],$$

which means that $K_\beta([a, b]) \subset [a, b]$. Similarly, we obtain that $U_\alpha([a, b]) \subset [a, b]$ and $V([a, b]) \subset [a, b]$.

As $T(x^*) = (VU_\alpha K_\beta)(x^*) = VU_\alpha(x^*) = V(x^*) = x^*$, we have

$$|T^n(x) - x^*| = |T^n(x) - T^n(x^*)| \tag{21}$$

$$\begin{aligned}
&= |(VU_\alpha K_\beta)^n(x) - (VU_\alpha K_\beta)^n(x^*)| \\
&\leq \gamma^n |(U_\alpha K_\beta)^n(x) - (U_\alpha K_\beta)^n(x^*)| \\
&\leq (\gamma\delta)^n |K_\beta^n(x) - K_\beta^n(x^*)| \\
&\leq (\gamma\delta\lambda)^n |x - x^*|.
\end{aligned}$$

Since $0 < \gamma\delta\lambda < 1$, on taking limit as $n \to \infty$ on the both sides of the above inequality, we have, $T^n(x_0) \to x^*$ for each $x_0 \in [a,b]$. Therefore, $T^N$ satisfies all the conditions of Lemma 1 and hence $x^*$ is the unique fixed point of $T^N$. Thus, $x^*$ is a fixed point of $T$.

(ii) From (12),

$$x_n - x_{n+1} = \beta \frac{f(x_n)}{f'(x_n)} + \frac{f(y_n)}{f'(y_n)} + \alpha \frac{f(z_n)}{f'(z_n)} \qquad (22)$$

We know that $x^*$ is the root of $f(x)$ from the proof of (i) in Theorem 3. By using the mean value theorem and (22),

$$x_n - x_{n+1} = \beta \frac{f'(a_n)(x_n - x^*)}{f'(x_n)} + \frac{f'(b_n)(y_n - x^*)}{f'(y_n)} + \alpha \frac{f'(c_n)(z_n - x^*)}{f'(z_n)} \qquad (23)$$

$$= \beta \frac{f'(a_n)(x_n - x^*)}{f'(x_n)} + \frac{f'(b_n)}{f'(y_n)}\left(1 - \alpha \frac{f'(c_n)}{f'(z_n)}\right)\left(1 - \beta \frac{f'(a_n)}{f'(x_n)}\right)(x_n - x^*)$$

$$+ \alpha \frac{f'(c_n)}{f'(z_n)}\left(1 - \beta \frac{f'(a_n)}{f'(x_n)}\right)(x_n - x^*)$$

$$= (x_n - x^*)\left[\beta \frac{f'(a_n)}{f'(x_n)} + \left(1 - \beta \frac{f'(a_n)}{f'(x_n)}\right)\left(\frac{f'(b_n)}{f'(y_n)}\left(1 - \alpha \frac{f'(c_n)}{f'(z_n)}\right) + \alpha \frac{f'(c_n)}{f'(z_n)}\right)\right]$$

where $a_n = x^* + \mu(x_n - x^*)$, $b_n = x^* + \mu(y_n - x^*)$ and $c_n = x^* + \mu(z_n - x^*)$, $0 < \mu < 1$.

From (23), we have

$$\frac{x_n - x^*}{x_n - x_{n+1}} = \frac{1}{\beta \frac{f'(a_n)}{f'(x_n)} + \left(1 - \beta \frac{f'(a_n)}{f'(x_n)}\right)\left[\frac{f'(b_n)}{f'(y_n)}\left(1 - \alpha \frac{f'(c_n)}{f'(z_n)}\right) + \alpha \frac{f'(c_n)}{f'(z_n)}\right]}. \qquad (24)$$

Using conditions $(f_3)$, we obtain that

$$\frac{x_n - x^*}{x_n - x_{n+1}} \leq \frac{1}{\beta\frac{m}{M} - \left[\frac{M}{m} + \alpha\frac{M}{m}\right]}$$

$$= \frac{mM}{(1+\alpha)M^2 - \beta m^2}.$$

Consequently, for $n \geq 0$,

$$|x_n - x^*| < \frac{mM}{(1+\alpha)M^2 - \beta m^2} |x_n - x_{n+1}|,$$

which is a required error estimation.

**Remark 2.** *(i) Note that the error estimate (14) is better than the error estimate (4) for $\beta < \alpha$ and $2m > M$. Indeed, for $2m > M$ and $\beta < \alpha$, since*

$$\frac{\frac{mM}{(1+\alpha)M^2 - \beta m^2}}{\frac{M}{m}} = \frac{m^2}{(1+\alpha)M^2 - \beta m^2}$$

$$< \frac{m^2}{(1+\beta)m^2 - \beta m^2} = 1$$

*we get*

$$\frac{mM}{(1+\alpha)M^2 - \beta m^2} < \frac{M}{m}.$$

*(ii) Also, we can see that the error estimate (14) is better than the error estimate (2.1) in [12] for $\beta < \alpha$ and $2m > M$. Indeed,*

$$\beta < \alpha \Rightarrow -\alpha < -\beta$$
$$\Rightarrow -\alpha m^2 < -\beta m^2$$
$$\Rightarrow M^2 - \alpha m^2 < M^2 - \beta m^2$$
$$\Rightarrow M^2 - \alpha m^2 < (1+\alpha)M^2 - \beta m^2$$

*we have*

$$\frac{mM}{(1+\alpha)M^2 - \beta m^2} < \frac{mM}{M^2 - \alpha m^2}.$$

## 3. Polynomiographs

Polynomiography bridges the gap between math and art, combining them into patterns that have symmetry and equilibrium. Polynomials themselves have wide uses in mathematics. Polynomiography extends those uses by allowing users to see clearly the basins of attraction and speed of convergence of a selected root-finding method. This can give greater insight into various classes of polynomials. Moreover, polynomiography has applications in a number of artistic practices, including design. Polynomiographs have been used as inspiration for many mediums, such as painting, sculpting and weaving.

Polynomials are undoubtedly one of the most significant objects in all of mathematics and sciences. The problem of polynomial roots finding was known since Sumerians 3000 years B.C. Over the centuries, mathemeticians have developed a variety of methods of solving equations. In 17th century Newton proposed a method for calculating approximately roots of polynomials. The behavior of Newton's method in the complex plane as applied to the equation $z^3 - 1 = 0$ investigated by Cayley in 1879 [6]. The Cayley's problem was solved by Julia in 1919 and then Mandelbrot in 1970 [17]. The last interesting contribution to the polynomials root finding was made by Kalantari [10]. Kalantari has developed visualization software that brings the process of finding the roots of a polynomial equation into the field of design and art. In 2005 he get U.S. patent for the technology of polynomiography [8].

Fractals and polynomiographs are obtained by iterations. Fractals are self-similar and independent of scale. This means there detail on all levels of magnification. On the other hand, polynomiography is well controlled and images of polynomiography are more predictable as compared to fractals. An infinite variety of designs can be created by using the infinite variety of complex polynomials.

According to the Fundamental Theorem of Algebra, any complex polynomial with comlex coefficients:

$$p(z) = a_n z^n + a_{n-1} z^{n-1} + \ldots + a_1 z + a_0 \tag{25}$$

of degree $n$ has $n$ roots. The degree $n$ of polynomial describes the number of basins of attraction in complex plane. Restuating the roots on the complex plane manually, localizations of basins can be controlled. Description of polynomiograph, its theoretical background and artistic applications are described in [9], [10].

In [15] Kotarski et al. used the Mann and Ishikawa iterations instead of the standart Picard iteration to obtain some generalization of Kalantari's polynomiography. They introduced some polynomiographs for the cubic equation $z^3 - 1 = 0$. Latif et al. in [19], using the ideas from [15], have used the S-iteration in polynomiography.

In this section we recall the well-known Newton method for finding roots of a complex polynomial $p$. The Newton method is given as followig:

$$z_{n+1} = N(z_n), n \geq 0. \tag{26}$$

where $N(z) = z - \frac{p(z)}{p'(z)}$ and $z_0 \in \mathbb{C}$ is a starting point. $p'(z)$ is the first derivative of $p$ at $z$. We will take the space $X = \mathbb{C}$ or $X = \mathbb{R}^2$ that is Banach one. We take $z_0 = (x_0, y_0)$ and $\alpha_n = \alpha$, $\beta_n = \beta$ such that $0 \leq \alpha \leq 1$ and $0 \leq \beta \leq 1$.

Applying the Picard-Mann hybrid iteration process (10) in (13) we obtain the following formula:

$$\begin{cases} z_{n+1} = N(v_n), \\ v_n = (1-\alpha)z_n + \alpha N(z_n), n \geq 0 \end{cases} \tag{27}$$

where $\alpha \in (0,1]$.

Using the iteration process (9) in (13) we get:

$$\begin{cases} z_{n+1} = (1-\alpha)v_n + \alpha N(v_n), \\ v_n = (1-\beta)w_n + \beta N(w_n), \\ w_n = (1-\gamma)z_n + \gamma N(z_n), n \geq 0 \end{cases} \tag{28}$$

where $\alpha \in (0,1]$ and $\beta, \gamma \in [0,1]$.

Substituting the our iteration (11) in (13) we get:

$$\begin{cases} z_{n+1} = N(v_n), \\ v_n = (1-\alpha)w_n + \alpha N(w_n), \\ w_n = (1-\beta)z_n + \beta N(z_n), n \geq 0 \end{cases} \qquad (29)$$

where $\alpha \in (0,1]$ and $\beta \in [0,1]$.

The sequence $\{z_n\}_{n=0}^{\infty}$ is called the orbit of the point $z_0$. If the sequence $\{z_n\}_{n=0}^{\infty}$ converges to a root $z^*$ then we say that $z_0$ is attracted to $z^*$. A set of all starting points for $\{z_n\}_{n=0}^{\infty}$ converges to $z^*$ is called the basin of attraction of $z^*$. Boundaries between basins usually are fractals in nature. The formulas given above are used in the next section to obtain polynomiographs for complex polynomials that visualize the roots finding process.

## 4. Examples of Polynomiographs with Different Iterations

In this section a few examples of the polynomiographs are obtained using iteration processes (27)-(29) defined in the previous sections are presented. In our experiments we focused on the comparison of the different iteration processes for discrete values of parameters. These polynomiographs for different parameters and different complex equations as follows:

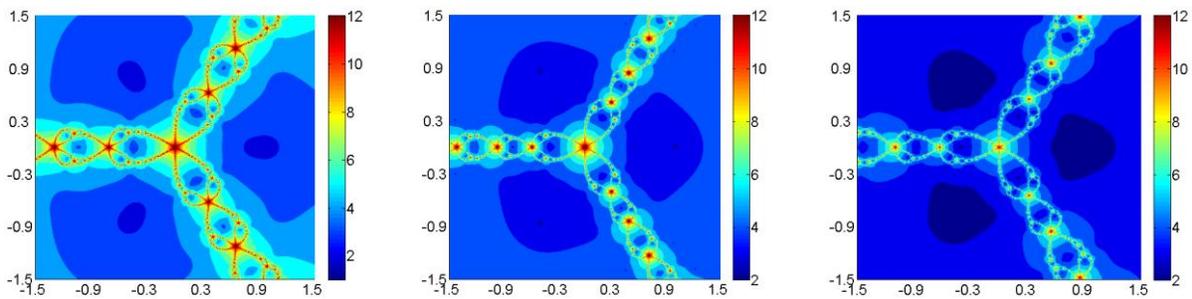

iteration (27)          iteration (28)          iteration (29)

**Figure 1:** Examples of polynomiographs for $z^3 - 1 = 0$

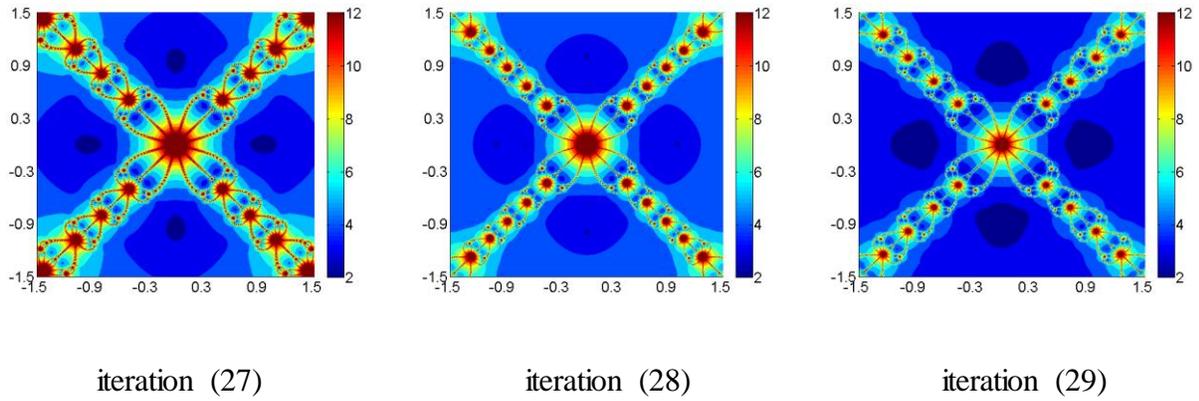

iteration (27)        iteration (28)        iteration (29)

**Figure 2:** Examples of polynomiographs for $z^4 - 1 = 0$

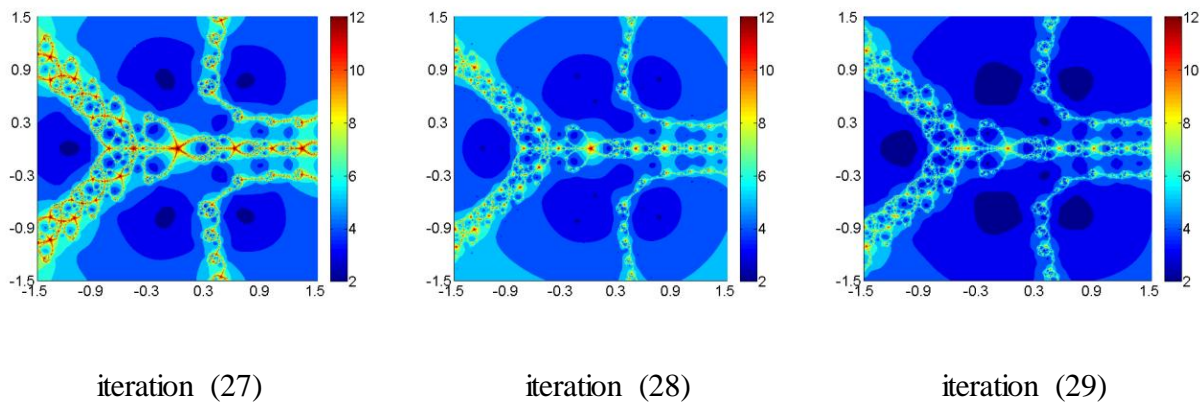

iteration (27)        iteration (28)        iteration (29)

**Figure 3:** Examples of polynomiographs for $z^5 + z^2 + 1 = 0$

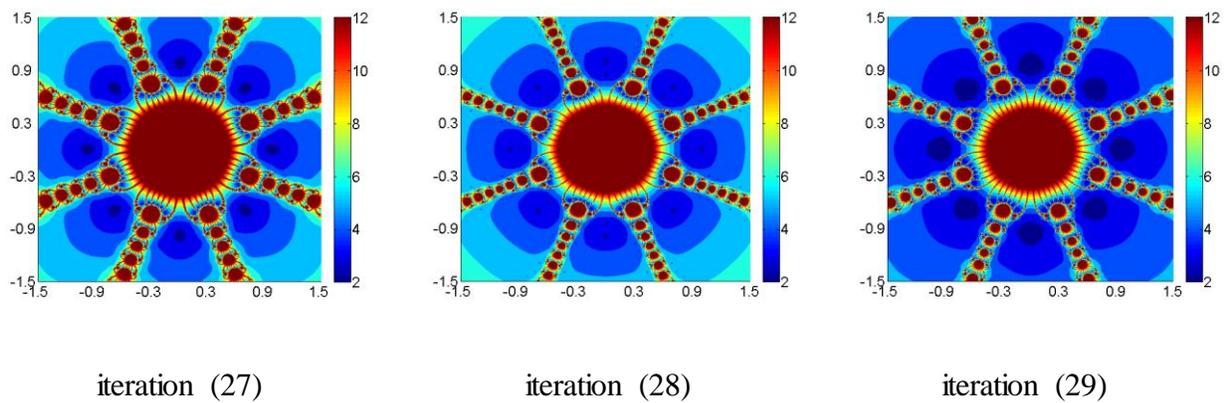

iteration (27)        iteration (28)        iteration (29)

**Figure 4:** Examples of polynomiographs for $z^8 - 1 = 0$

In Figure 1 three images with three basins of attraction to the three roots of polynomial $z^3 - 1 = 0$, in Figure 2 three images with four basins of attraction to the four roots of polynomial $z^4 - 1 = 0$, in Figure 3 three images with five basins of attraction to the five roots of polynomial $z^5 + z^2 + 1 = 0$, in Figure 4 three images with eight basins of attraction to the eight roots of polynomial $z^8 - 1 = 0$ are presented. These equations were solved in the square $[-1.5, 1.5] \times [-1.5, 1.5]$ using three different iteration processes described in the previous section. The different colours of image depend on the number of iterations needed to reach a root with the given accuracy $\varepsilon = 0.001$. The upper bound of the number of iterations was fixed as $k = 12$ and the parameters used in the iterations were fixed as $\alpha = 0.8$, $\beta = 0.6$ and $\gamma = 0.5$. By changing parameters $\alpha$, $\beta$, $\gamma$, $\varepsilon$ and $k$ one can obtain infinitely many polynomiographs.

All the experiments were performed on a computer with the following specification: Intel Core i3 processor, 2.53 GHz, 4GB RAM and Windows 7 (64-bit). MATLAB software was used for generating polynomiographs.

**Acknowledgement 1.** *This work was supported by Ataturk University Rectorship under " The Scientific and Research Project of Ataturk University" , Project No.: 2016/153.*

**References**

[1] Agarwal, R.P., O'Regan, D., Sahu, D.R., Iterative construction of fixed points of nearly asymptotically nonexpansive mappings, J. Nonlinear Convex Anal. 8 (2007) 61–79.

[2] Berinde, V., Conditions for the convergence of the Newton method, An. St. Univ. Ovidius Constanta, **3**(1995), No. 1, 22-28.

[3] Berinde, V., On some exit criteria for the Newton method, Novi Sad J. Math., **27**(1997), No. 1, 19-26.

[4] Berinde, V., On the extended Newton's method, in Advances in Difference Equations, S. Elaydi, I. Gyori, G. Ladas (eds.), Gordon and Breach Publishers, 1997, 81-88.

**[5]** Berinde, V., Păcurar, M., A fixed point proof of the convergence of a Newton-type method, Fixed Point Theory, **7** (2006), No. 2, 235-244

**[6]** Cayley, A., The Newton-Fourier Imaginary Problem, American Journal of Mathematics, 2 (1879), p. 97.

**[7]** Demidovich, B. P., Maron, A. I., Computational Mathematics, MIR Publishers, Moscow, 1987.

**[8]** Kalantari, B., Method of Creating Graphical Works Based on Polynomials, U.S. Patent 6,894,705, (2005).

**[9]** Kalantari, B., Polynomiography: From the Fundamental theorem of Algebra to Art, Leonardo, 38 (2005), No. 3, 233-238.

**[10]** Kalantari, B., Polynomial Root-Finding and Polynomiography, World Scientific, Singapore (2009).

**[11]** Kantorovich, L. V., Akilov, G. P., Functional analysis, Second edition, Pergamon Press, Oxford-Elmsford, New York, 1982.

**[12]** Karaca, N., Abbas, M. and Yildirim, I., Convergence of a Newton-Like S-Iteration Process in $\mathbb{R}$, Creative Math. Inf., (In Press).

**[13]** Karaca, N., Yildirim, I., Approximating Fixed Points of Nonexpansive Mappings by a Faster Iteration Process, J. Adv. Math. Stud., Vol. 8 (2015), No. 2, 257-264.

**[14]** Khan, S. H., A Picard-Mann hybrid iterative process, Fixed Point Theory and Applications, vol. 2013, article 69, 10 pages, 2013.

**[15]** Kotarski, W., Gdawiec, K., Lisowska, A., Polynomiography via Ishikawa and Mann Iterations, G. Bebis et al. (eds.) Advances in Visual Computing, Part I, LNCS 7431, Springer, Berlin, (2012) 305-313.

**[16]** Latif, A., Rafiq, A., Shahid, A. A., Polynomiography via S-iteration Scheme, Abstract and Applied Analysis, In Press.

**[17]** Mandelbrot, B., The Fractal Geometry of Nature. W.H. Freeman and Company, New York (1983).

**[18]** Ortega, J., Rheinboldt, W. C., Iterative solution of nonlinear equations in several variables, Academic Press, New York, 1970.

**[19]** Scherzer, O., Convergence criteria of iterative methos based on Landweber iteration for solving nonlinear problems, J. Math Anal Appl. 194, 911–933 (1995). doi:10.1006/jmaa.1995.1335

**[20]** Sen, R. N., Biswas, A., Patra, R., Mukherjee, S., An extension on Berinde's criterion for the convergence of a Newton-like method, Bull. Calcutta Math. Soc. (to appear).

**[21]** Sen, R. N., Mukherjee, S., Patra, R., On the convergence of a Newton-like method in R and the use of Berinde's exit criterion, Intern. J. Math. Math. Sc. Vol. 2006 (2006), Article ID 36482, 9 pages; doi:10.1155/IJMMS/2006/36482.

**[22]** Thianwan, S., Common fixed points of new iterations for two asymptotically nonexpansive nonself mappings in a Banach space, J. Comput. Appl. Math., In Press, (2008), doi:10.1016/j.cam.2008.05.051.

**[23]** Yildirim, I., Ozdemir, M., A new iterative process for common fixed points of finite families of non-self-asymptotically non-expansive mappings, Nonlinear Analysis, 71 (2009) 991-999.

**DEPARTMENT OF MATHEMATICS, FACULTY OF SCIENCE, ATATURK UNIVERSITY, 25240 ERZURUM, TURKEY**

**E-mail address : nazli.kadioglu@atauni.edu.tr**

**E-mail address : isayildirim@atauni.edu.tr**